\documentclass{article}
\usepackage[english]{babel}
\usepackage{amsmath}
\usepackage{amsthm}
\usepackage{amssymb}
\usepackage{amscd}
\usepackage[latin1]{inputenc}

\newtheorem{teor}{Theorem}
\newtheorem{defi}{Definition}
\newtheorem{lema}{Lemma}

\newtheorem{cor}{Corollary}
\newtheorem{rem}{Remark}
\newtheorem{ejem}{Example}

\begin{document}

\title{On semiprimary rings of finite global dimension}

\author{Manuel Saorín \\ Departamento de Matemáticas \\
Universidad de Murcia, Aptdo. 4021\\ 30100 Espinardo, Murcia\\
SPAIN\\ {\it e-mail}: msaorinc@um.es}

\date{}

\thanks{We warmly thank Dan Zacharia for very fruitful
conversations and for telling us about the problem tackled in the
paper. We also thank the D.G.I. of the Spanish Ministry of Science
and  Technology and the fundación 'Séneca' of Murcia for their
financial support}

\maketitle

\begin{abstract}

{\bf  Suppose that $A$ is a semiprimary ring satisfying one of the
two conditions: 1) its Yoneda ring is generated in finite degrees;
2) its Loewy length is less or equal than three. We prove that the
global dimension of $A$ is finite if, and only if, there is a
$m>0$ such that $Ext_A^n(S,S)=0$, for all simple $A$-modules $S$
and all $n\geq m$}

\end{abstract}

In a recent paper, Skowronski, Smal$\emptyset$ and Zacharia
(\cite{SSZ}) proved that a left Artinian ring $A$ has finite
global dimension if, and only if, every finitely generated
indecomposable left $A$-module has either finite projective
dimension or finite injective dimension. Then, they  showed that
one cannot replace 'indecomposable' by 'simple' in that statement,
by giving a counterexample of Loewy length 4. Finally they asked
the following question:

{\it Question}: Suppose $A$ is a left Artinian ring such that, for
each finitely generated indecomposable left $A$-module $M$, one
has $Ext_A^n(M,M)=0$ for $n\gg 0$ (i.e. there is a $m=m(M)>0$ such
that $Ext_A^n(M,M)=0$ for all $n\geq m$). Does $A$ have finite
global dimension?

The two main results of this paper, Theorems \ref{finite-degrees}
and \ref{radical-cube},  give two partial affirmative answers to
the latter question, even in the more general context of
semiprimary rings, showing in the way that the above mentioned
counterexample is of minimal Loewy length.

All rings in the paper are associative with unit. A ring $A$ is
{\bf semiprimary} when its Jacobson radical $J=J(A)$ is nilpotent
and  $A/J$ is semisimple. In that case, the minimal $n$ such that
$J^n=0$ is called the {\bf Loewy length} of $A$. Recall that a
semiprimary ring is a particular instance of (left and right)
perfect ring and, hence, every $A$-module $M$ has a projective
cover $\epsilon :P_0=P(M)\longrightarrow M$ and a minimal
projective resolution
$...P_{n+1}\stackrel{d_{n+1}}{\longrightarrow}P_n
\stackrel{d_n}{\longrightarrow}...P_1\stackrel{d_1}{\longrightarrow}P_0
\stackrel{\epsilon}{\longrightarrow}{M}\rightarrow 0$, both
uniquely determined up to isomorphism. We shall put
$\Omega^nM=Im(d_n)=Ker(d_{n-1})$, for all $n>0$, and
$\Omega^0M=M$. By dimension shifting, we have $Ext_A^n(M,X)\cong
Ext_A^1(\Omega^{n-1}M,X)$, for every $A$-module $X$ and every
$n>0$. When $X$ is semisimple and we apply the functor
$Hom_A(-,X)$ to the canonical exact sequence
$0\rightarrow\Omega^nM\stackrel{j}{\hookrightarrow}
P_{n-1}\longrightarrow\Omega^{n-1}M\rightarrow 0$, the induced map
$j^*:Hom_A(P_{n-1},X)\longrightarrow Hom_A(\Omega^nM,X)$ is the
zero map, so that $Ext_A^1(\Omega^{n-1}M,X)\cong
Hom_A(\Omega^nM,X)$, for all $n>0$. From that we get that
$Ext_A^n(M,X)\cong Hom_A(\Omega^nM,X)$, for all $n\geq 0$.  All
throughout the paper, whenever necessary, we shall see the latter
isomorphism as an identification. On the other hand, the
perfectness of $A$ implies that flat and projective modules
coincide, from which it follows that the left and right global
dimensions coincide with the weak global dimension, which is
left-right symmetric (cf \cite{Ro}[Theorem 9.15]). We shall deal
only with left modules, hitherto just called 'modules', although
the results obtained are left-right symmetric. Since the radical
filtration on any module is finite and has semisimple factors, the
global dimension of $A$ is the supremum of the projective
dimensions of  the simple $A$-modules. As a consequence, that
global dimension is finite if, and only if, $Ext_A^n(A/J,A/J)=0$
for some $n>0$.

Recall that, for $A$-modules $M,N,P$ and $m,n\geq 0$, one has a
$\mathbf{Z}$-bilinear {\bf Yoneda product} $Ext_A^m(N,P)\times
Ext_A^n(M,N)\longrightarrow Ext_A^{m+n}(M,P)$ (cf.
\cite{ML}[Chapter III, Section 5]). We shall need an explicit
description of this product when $N,P$ are semisimple,  using the
above mentioned identification. In this case, if $f\in
Hom_A(\Omega^nM,N)\cong Ext_A^n(M,N)$, then the comparison theorem
(cf. \cite{ML}[Theorem III.6.1]) yields a chain map between the
minimal projective resolutions of $\Omega^nM$ and $M$, which is
uniquely determined up to homotopy. Then we get a morphism
$\tilde{f}:\Omega^{m+n}M\longrightarrow\Omega^mN$. If now $g\in
Hom_A(\Omega^mN,P)\cong Ext_A^m(N,P)$ then the Yoneda product
$g\cdot f$ is just the composition $g\circ\tilde{f}$ and does not
depend on the choices made to get $\tilde{f}$. It is well-known
that the Yoneda product makes $E=Ext_A^*(A/J,A/J)=\oplus_{n\geq
0}Ext_A^n(A/J,A/J)$ into a graded ring,  with the obvious grading.
It is called the {\bf Yoneda ring} of $A$ and
$Ext_A^*(M,A/J)=\oplus_{n\geq 0}Ext_A^n(M,A/J)$ is canonically a
graded (left) $E$-module, for every $A$-module $M$. In case $A$ is
an algebra over a commutative ring $K$, the grading of $E$ makes
it into a graded $K$-algebra.

\vspace*{0.3cm}
\begin{defi}
A positively graded ring $R=\oplus_{n\geq 0}R_n$ will be said to
be {\bf generated in finite degrees} when there is a $m\geq 0$
such that the subgroup $R_0\oplus ...\oplus R_m$ generates $R$ as
a ring, i.e.,  there is no proper subring of $R$ containing
$R_0\oplus ...\oplus R_m$
\end{defi}

\begin{rem}
If $K$ is a commutative ring and $R=\oplus_{n\geq 0}R_n$ is a
positively graded $K$ algebra such that each $R_n$ is a finitely
generated $K$-module, then $R$ is generated in finite degrees if,
and only if, $R$ is finitely generated as a $K$-algebra
\end{rem}

The following is the first main result of the paper:

\begin{teor}
\label{finite-degrees} Let  $A$ be a semiprimary ring. The
following assertions are equivalent:

\begin{enumerate}
\item The global dimension of $A$ is finite \item The Yoneda ring
$E=Ext^*(A/J,A/J)$ is generated in finite degrees and
$Ext^n_A(S,S)=0$, for all simple  $A$-modules $S$ and all $n>>0$
\end{enumerate}
\end{teor}

{\bf Proof:} We only need to prove $2)\Longrightarrow 1)$. Let us
fix some $s>0$ such that $\oplus_{0\leq i\leq s}Ext_A^i(A/J,A/J)$
generates $E$ as a ring. We denote by $r$ the number of
nonisomorphic simple $A$-modules and fix $m>Sup\{k\geq 0:$
$Ext_A^k(T,T)\neq 0,$ for some simple $_AT\}$, something which is
possible to do since our hypothesis guarantees the existence of
such a supremum. We claim that if $n>m\cdot r\cdot s$ then
$Ext_A^n(A/J,A/J)=0$ and, hence,  assertion 1) will follow.
Indeed, if $n>m\cdot r\cdot s$ then $Ext_A^n(A/J,A/J)$ is
contained in a sum of products of the form
$Ext_A^{i_1}(A/J,A/J)\cdot ...\cdot Ext_A^{i_t}(A/J,A/J)$, with
$i_1+...+i_t=n$ and  $i_k\leq s$ for $k=1,...,t$. From that we get
that $n\leq t\cdot s$, which implies
$t\geq\frac{n}{s}>\frac{m\cdot r\cdot s}{s}=m\cdot r$. But
$Ext_A^{i_1}(A/J,A/J)\cdot ...\cdot Ext_A^{i_t}(A/J,A/J)$ is
contained in a sum of products of the form
$Ext_A^{i_1}(S_{t-1},S_t)\cdot ...\cdot Ext_A^{i_t}(S_0,S_1)$,
with the $S_i$ simple. In any such product, there exists a simple
$S_i$ which appears repeated, at least, $m+1$ times. Suppose
$j_0<j_1<...<j_m$ are different indices of $\{0,1,...,t\}$ such
that $S_{j_0}=S_{j_1}=...=S_{j_m}$, a simple $A$-module which we
denote by $X$. In case the product $Ext_A^{i_1}(S_{t-1},S_t)\cdot
...\cdot Ext_A^{i_t}(S_0,S_1)$ is nonzero, one gets that
$Ext_A^j(X,X)\neq 0$, for some $j\geq m$, which is a
contradiction. As a consequence, $Ext_A^n(A/J,A/J)=0$ as desired.

\vspace*{0.3cm} Using now the remark, the following consequence is
straightforward.

\begin{cor}
Let $A$ be an Artin algebra with center $K$. The following
assertions are equivalent:

\begin{enumerate}
\item The global dimension of $A$ is finite \item
$E=Ext_A^*(A/J,A/J)$ is finitely generated as a $K$-algebra and
$Ext_A^n(S,S)=0$, for all simples  $A$-modules $S$ and all $n\gg
0$
\end{enumerate}
\end{cor}

\begin{ejem}
Let $K$ be a field and $A$ the monomial algebra with quiver $Q:
1\rightleftarrows 2$ and the path
$1\stackrel{\alpha}{\longrightarrow}2\stackrel{\beta}{\longrightarrow}1
\stackrel{\alpha}{\longrightarrow}2$ as unique zero-relation. By
\cite{SSZ}, we know that $gl.dim(A)=\propto$ but $Ext_A^n(S,S)=0$
for all simples $S$ and all $n\gg 0$. According to the above
theorem, the Yoneda algebra $E=Ext_A^*(A/J,A/J)$ cannot be
finitely generated. That can be explicitly seen by using Anick and
Green's resolution (cf. \cite{AG}). Using the terminology of this
latter paper, the sets of chains are $\Gamma_n=Q_n$ ($n=0,1$) and
$\Gamma_n=\{\alpha (\beta\alpha)^{n-1}\}$, for all $n\geq 2$. By
using \cite{GZ}[Theorem A, Proposition 1.2], we see that $E$ can
be represented by a quiver $\tilde{Q}$ having $\tilde{Q}_0=Q_0$,
one arrow $b:1\rightarrow 2$ corresponding to the 1-chain $\beta$
and arrows $a_n:2\rightarrow 1$ ($n>0$), where $a_n$ corresponds
to the $n$-chain $\alpha (\beta\alpha)^{n-1}$. All possible
products of arrows in $\tilde{Q}$ are zero and the grading on $E$
assigns degree zero to the vertices,  $deg(b)=1$ and $deg(a_n)=n$,
for all $n>0$.
\end{ejem}

\newpage

We move now to study the case of Loewy length $\leq 3$. The
following lemma is probably known, but we include a proof for the
sake of completeness.

\begin{lema}
\label{auxiliar} Let $A$ be a ring such that $A/J$ is semisimple
and $J^3=0$, where $J=J(A)$. Let $_AP$ be a projective module and
$M$ a submodule of $JP$. Then $M$ decomposes as $M=N\oplus X\oplus
Y$, where $X$ is a semisimple direct summand of $JP$, Y is a
(semisimple) direct summand of $J^2P$ and $N$ is a submodule
containing no simple direct summand. Moreover, the following
assertions are equivalent:

\begin{enumerate}
\item $M\cap J^2P=JM$ \item $Y=0$ in the above decomposition
\end{enumerate}
\end{lema}

{\bf Proof:} Since $J^2M=0$, $JM$ is a semisimple submodule of $M$
and we can fix a decomposition $Soc(M)=JM\oplus Z$. Then the
canonical composition
$Z\stackrel{j}{\hookrightarrow}M\stackrel{p}{\longrightarrow}M/JM$
is a (necessarily split) mono. If now $\pi :M/JM\longrightarrow Z$
is a retraction for $p\circ j$, then $\pi\circ p$ is a retraction
for $j$, so that, putting $N=Ker(\pi\circ p)$, we get a
decomposition $M=N\oplus Z$. We claim that $N$ contains no simple
direct summand. To see that, notice first that $JM=JN$. If, by of
contradiction, we assume that there exists a simple direct summand
$S$ of $N$, then we have a decomposition $N=S\oplus N'$, from
which follows that $JM=JN=JN'$. Now we get a decomposition
$Soc(M)=Soc(N')\oplus S\oplus Z$, with $JM\subseteq Soc(N')$. That
contradicts the initial decomposition. Now, since $Z$ is
semisimple, we have a decomposition $Z=(Z\cap J^2P)\oplus X$. But
$X\cap J^2P=0$, which implies that the canonical composition
$X\hookrightarrow JP\stackrel{pr}{\longrightarrow}JP/J^2P$ is a
(split) mono and, arguing as above, we conclude that $X$ is a
(semisimple) direct summand of $JP$. Next we prove the equivalence
of assertions 1 and 2:

$1)\Longrightarrow 2)$ We always have inclusions $JM=JN\subseteq
JN\oplus Y\subseteq M\cap J^2P$. Hence, assertion 1) implies that
$Y=0$

$2)\Longrightarrow 1)$ We just need to prove the inclusion $M\cap
J^2P\subseteq JM$. Take $m\in M\cap J^2P$, which we decompose as
$m=n+x$, with $n\in N$ and $x\in X$. Since $Jm=0=Jx$ we get that
$Jn=0$ and $An$ is a semisimple submodule of $N$. We claim that
$n\in JM$. Indeed, if $n\notin JM=JN$, then the semisimplicity of
$An$ yields a decomposition $An=(An\cap JN)\oplus Aan$, with
$an\neq 0$. Now the canonical composition $Aan\hookrightarrow
N\stackrel{p}{\longrightarrow}N/JN$ is monic and, arguing as
above, we get that $Aan$ is a nonzero (semisimple) direct summand
of $N$, which is a contradiction, thus settling our claim. Since
$n\in JM\subseteq M\cap J^2P$, we get $x\in J^2P$. But, since $Ax$
is a semisimple direct summand of $JP$, the latter implies that
$x=0$ and, hence, $m=n\in JM$.

\vspace*{0.5cm} Using the above lemma, to each simple $A$-module
$S$ we can inductively associate a sequence of decompositions
$\Omega^nS=M_n\oplus Z_n$ as follows. We put $M_n=\Omega^nS$,
$Z_n=0$ for $n=0,1$. If $n>1$ and the decomposition
$\Omega^{n-1}S=M_{n-1}\oplus Z_{n-1}$ is already defined, then we
put $P'_{n-1}=P(M_{n-1})$ and $P''_{n-1}=P(Z_{n-1})$. According to
the above lemma, we have a decomposition $\Omega M_{n-1}=M_n\oplus
Y_n$, where $Y_n$ is a direct summand of $J^2P'_{n-1}$ and
$M_n\cap J^2P'_{n-1}=JM_n$. Then, by putting $Z_n=Y_n\oplus\Omega
Z_{n-1}$, we get the desired decomposition $\Omega^nS=M_n\oplus
Z_n$. In this way we get a sequence of modules
$(S=M_0,M_1,M_2,...)$ which is uniquely determined by $S$. When we
look at the canonical isomorphism $Ext_A^n(S,A/J)\cong
Hom_A(\Omega^nS,A/J)$ as an identification, we can view
$Hom_A(M_n,A/J)$ as a subgroup of $Ext_A^n(S,A/J)$. That is the
sense of the following crucial result.

\begin{lema}
\label{1-extensions} Let $A$ be a ring such that $A/J$ is
semisimple and $J^3=0$. With the above notation, $Hom_A(M_n,A/J)$
is contained in $Ext_A^1(A/J,A/J)\cdot\stackrel{n}{....}\cdot
Ext_A^1(A/J,A/J)\cdot Ext_A^0(S,A/J)$
\end{lema}

{\bf Proof:} Our argument is inspired by that of
\cite{GMV}[Proposition 3.2]. We  apply induction on $n$, the case
$n=0$ being trivially true. Suppose $n>0$. We plan to prove that
$Hom_A(M_n,A/J)\subseteq Ext_A^1(A/J,A/J)\cdot Hom_A(M_{n-1},A/J)$
(after the above mentioned identifications) and the induction
hypothesis will give the desired result. Let $f\in
Hom_A(M_n,A/J)$. Then $(f$ $0):M_n\oplus
Z_n=\Omega^nS\longrightarrow A/J$ represents an element of
$Ext_A^n(S,A/J)$. The canonical inclusion
$\iota:\Omega^nS=M_n\oplus Z_n\hookrightarrow JP'_{n-1}\oplus
JP''_{n-1}$ can be written as a matrix $\iota$ =\[ \left(
\begin{array}{lr}
i_{11} & i_{12}\\
0 & i_{22}
\end{array} \right) \]
where $i_{11}:M_n\hookrightarrow JP'_{n-1}$ is the inclusion,
$i_{22}$ is the map $(0$ $j):Z_n=Y_n\oplus\Omega
Z_{n-1}\longrightarrow JP''_{n-1}$, with $j:\Omega
Z_{n-1}\longrightarrow JP''_{n-1}$ the canonical inclusion, and
$i_{12}$ is the map $(i$ $0): Z_n=Y_n\oplus\Omega
Z_{n-1}\longrightarrow JP'_{n-1}$, with $i:Y_n\hookrightarrow
JP'_{n-1}$ the canonical inclusion. Notice that, by definition, we
have  $Im(i_{12})\subseteq J^2P'_{n-1}$. As a consequence, the
induced homomorphism $\Omega^nS/J\Omega^nS=(M_n/JM_n)\oplus
(Z_n/JZ_n)\longrightarrow (JP'_{n-1}/J^2P'_{n-1})\oplus
(JP''_{n-1}/J^2P''_{n-1})$ has a diagonal matrix shape \[ \left(
\begin{array}{lr}
\bar{i}_{11} & 0\\
0 & \bar{i}_{22}
\end{array} \right) \].

 On the other hand, $f$ factors as $f=\bar{f}\circ
p_M$, where $p_M:M_n\longrightarrow M_n/JM_n$ is the canonical
projection, and $\bar{i}_{11}:M_n/JM_n\longrightarrow
JP'_{n-1}/J^2P'_{n-1}$ is a (split) monomorphism by the definition
of $M_n$. If we choose a retraction $\alpha$ for $\bar{i}_{11}$
and put $\bar{g}=\bar{f}\circ\alpha$, then the induced morphism
$(\bar{g}$ $0):(JP'_{n-1}/J^2P'_{n-1})\oplus
(JP''_{n-1}/J^2P''_{n-1})\longrightarrow A/J$ has the property
that $(\bar{g}$ $0)\circ$ \[ \begin{pmatrix} \bar{i}_{11} & 0\\
0 & \bar{i}_{22}
\end{pmatrix}\] $=(\bar{f}$ $0)$.

 We now denote by $g$ the composition
$JP'_{n-1}\stackrel{pr}{\longrightarrow}JP'_{n-1}/J^2P'_{n-1}
\stackrel{\bar{g}}{\longrightarrow}A/J$ and then the morphism $(g$
$0):JP'_{n-1}\oplus JP''_{n-1}\longrightarrow A/J$ has the
property that $(g$ $0)\circ\iota =(f$ $0)$. In particular, $g\circ
i_{11}=f$ and $g\circ i_{12}=0$. By consicering $(i_{11}$
$i_{12})$ with arrival in $P'_{n-1}$ instead of $JP'_{n-1}$,  we
have $M_{n-1}\cong Coker[(i_{11}$ $i_{12})]$ and, since the
composition $\pi\circ (i_{11}$ $i_{12})$ is zero, where $\pi
:P'_{n-1}\longrightarrow P'_{n-1}/JP'_{n-1}$ is the canonical
projection, we have a unique morphism $h:M_{n-1}\longrightarrow
P'_{n-1}/JP'_{n-1}$ such that the composition
$P'_{n-1}\stackrel{pr}{\longrightarrow}
M_{n-1}\stackrel{h}{\longrightarrow}P'_{n-1}/JP'_{n-1}$ is the
canonical projection. Now we have a commutative diagram:

\setlength{\unitlength}{1mm}
\begin{picture}(85,45)
\put(0,39){0} \put(2,40){\vector(1,0){7}} \put(10,39){$\Omega^nS$}
\put(17,40){\vector(1,0){14}} \put(32,39){$P'_{n-1}\oplus
P''_{n-1}$} \put(50,40){\vector(1,0){14}}
\put(65,39){$\Omega^{n-1}S$} \put(76,40){\vector(1,0){8}}
\put(85,39){0}

\put(0,19){0}\put(2,20){\vector(1,0){7}} \put(10,19){$JP'_{n-1}$}
\put(19,20){\vector(1,0){18}} \put(38,19){$P'_{n-1}$}
\put(47,20){\vector(1,0){18}}
\put(66,19){$\frac{P'_{n-1}}{JP'_{n-1}}$}
\put(76,20){\vector(1,0){8}} \put(85,19){$0$}

\put(10,0){$A/J$} \put(14,38){\vector(0,-1){15}}
\put(42,38){\vector(0,-1){15}} \put(71,38){\vector(0,-1){15}}
\put(14,18){\vector(0,-1){14}}

\put(15,29){$(i_{11}$ $i_{12})$} \put(15,9){$g$} \put(43,29){$(1,$
$0)$} \put(72,29){$(h,$ $0)$}

\end{picture}

\vspace*{0.3cm} where the rows are  the obvious exact sequences
and the composition of the vertical left arrows is $(f,$ $0)$.
Now, since $JP'_{n-1}=\Omega(P'_{n-1}/JP'_{n-1})$, we have $g\in
Ext_A^1(P'_{n-1}/JP'_{n-1},A/J)$  and, on the other hand, we also
have  $(h$ $0)\in Ext_A^{n-1}(S,P'_{n-1}/JP'_{n-1})$. Then the
Yoneda product $g\cdot (h$ $0)$ is just $g\circ (i_{11}$
$i_{12})=(f$ $0)$. Moreover, $(h$ $0)$ clearly belongs to
$Hom_A(M_{n-1},P'_{n-1}/JP'_{n-1})$ when we view the latter as a
subgroup of $Ext_A^{n-1}(S,P'_{n-1}/JP'_{n-1})$. From that we get
the desired inclusion $Hom_A(M_n,A/J)\subseteq
Ext_A^1(A/J,A/J)\cdot Hom_A(M_{n-1},A/J)$.

\vspace*{0.5cm} We are now in a position to prove:

\begin{teor}
\label{radical-cube} Let $A$ be a ring such that $A/J$ is
semisimple and $J^3=0$. The following assertions are equivalent:
\begin{enumerate}
  \item The  global dimension of $A$ is finite
\item $Ext_A^n(S,S)=0$, for all simple  $A$-modules $S$ and all
$n\gg 0$ \item Every simple $A$-module has either finite
projective dimension or finite injective dimension
\end{enumerate}
\end{teor}

{\bf Proof:} We only need to prove $2)\Longrightarrow 1)$. If
condition 2) holds, even without the hypothesis that $J^3=0$, one
can define an order relation among the simple modules by the rule:
$S\preceq T$ iff $S$ is a direct summand of $\Omega^mT$, for some
$m\geq 0$. Indeed, that is always a preorder relation and we only
need to check antisymmetry. If $S\preceq T\preceq S$ then $S$ is a
direct summand of $\Omega^mT$ and $T$ is a direct summand of
$\Omega^nS$, for some $m,n\in\mathbf{N}$. From that it follows
that $S$ is a direct summand of $\Omega^{m+n}S$, and hence a
direct summand of $\Omega^{(m+n)k}S$ for all $k\in\mathbf{N}$. If
$m+n>0$ we get a contradiction with the fact that $Ext_A^r(S,S)=0$
for all $r\gg 0$. Therefore $m+n=0$ or, equivalently, $m=n=0$,
which means that $S=T$ as desired.

We next claim that, for every simple module $S$, the sequence
$(S=M_0,M_1,M_2,...)$ defined immediately before the previous
lemma is finite, i.e., there is a $t>0$ such that  $M_t=0$ (and,
hence, $M_n=0$ for all $n\geq t$).  To see that, we follow an
argument analogous to that of Theorem \ref{finite-degrees}. Let
$r$  be the number of nonisomorphic simple $A$-modules and fix
$m>Sup\{k\geq 0:$ $Ext_A^k(T,T)\neq 0,$ for some simple $_AT\}$.
Now we claim that if $n>m\cdot r$, then $M_n=0$. Indeed, if
$M_n\neq 0$ then $Hom_A(M_n,A/J)\neq 0$ and the foregoing lemma
implies that $Ext_A^1(A/J,A/J)\cdot\stackrel{n}{....}\cdot
Ext_A^1(A/J,A/J)\neq 0$, from which we deduce the existence  of
simple $A$-modules $S_0,S_1,...,S_n$ such that
$Ext_A^1(S_{n-1},S_n)\cdot ...\cdot Ext_A^1(S_0,S_1)\neq 0$. But
the sequence of simples $(S_0,S_1,...,S_n)$ will necessarily have
one simple
 repeated, at least, $m+1$ times and, as in the proof of Theorem
 \ref{finite-degrees}, we get a simple $X$ such that
 $Ext_A^{j}(X,X)\neq 0$ for some $j\geq m$, which is a contradiction.
  Observe that,  when $S$ is minimal with respect
 to the order relation $\preceq$, one has $\Omega^nS=M_n$ for all
 $n\geq 0$. Hence, our argument also proves that $pd_A(S)<\propto$
 in that case.

 Let now $S$ be an arbitrary simple module and let us recall
 that, by definition of the associated sequence $(S=M_0,M_1,M_2,...)$,
 we have $\Omega M_{n-1}=M_n\oplus Y_n$, where $Y_n$ is a direct
 summand $J^2P(M_{n-1})$, for all $n>0$. From that we
 get, for every fixed $n>0$,  that
 $pd_A(S)<\propto$ if, and only if, the dimensions $pd(M_n)$,
 $pd_A(Y_1)$,...,$pd_A(Y_n)$ are all finite. Bearing in mind  the
 finiteness of the sequence $(M_n)$, we conclude that
 $pd_A(S)\propto$ iff $pd_A(Y_n)\propto$, for all $n>0$. But each
 $Y_n$ is a direct sum of simple modules which are strict
 predecessors of $S$ with respect to the order relation $\preceq$.
 Since the minimal
 elements with respect to $\preceq$ have finite projective
 dimension, an easy induction procedure shows that
 $pd_A(S)<\propto$, for every simple left $A$-module and, hence,
 $gl.dim (A)<\propto$.

\end{document}